\documentclass[12pt]{amsart}  
 \usepackage{amsmath}
 \usepackage{amssymb}
 \usepackage{amscd} 
 \usepackage[mathscr]{eucal}
 
 \newtheorem{thm}{Theorem}[section]
 
 \newtheorem{remark}[thm]{Remark}
 \newtheorem{example}[thm]{Example}   
 \newtheorem{cor}[thm]{Corollary}
 \newtheorem{prop}[thm]{Proposition}
 
 \newtheorem{defn}[thm]{Definition}

\newcommand{\R}{\mathbb{R}} 
 \newcommand{\Z}{\mathbb{Z}}

 \newcommand{\N}{\mathbb{N}}

 \renewcommand{\S}{\mathscr{S}}
 
 \newcommand{\A}{\mathcal{A}}
 \newcommand{\B}{\mathcal{B}}
 \newcommand{\V}{\mathcal{V}}
 \newcommand{\K}{\mathbf{k}}
 
 \newcommand{\tM}{\tilde{M}} 
 \newcommand{\lcm}{\operatorname{lcm}}

 \title{Topological representations of matroids} 
 \author{E. Swartz} 
 \address{Malott Hall \\ Cornell University \\ Ithaca, NY 14853} 
 \email{ebs@math.cornell.edu}
 \thanks{Partially supported by a VIGRE postdoc under NSF grant number 9983660
to Cornell University.}

\begin{document}
 \begin{abstract}
    There is a one-to-one correspondence between geometric lattices and the intersection lattices of arrangements of homotopy spheres.  When the arrangements are essential and fully partitioned, Zaslavsky's enumeration of the cells of the arrangement still holds. Bounded subcomplexes of an arrangement of homotopy spheres correspond to minimal cellular resolutions of the dual matroid Steiner ideal.  As a result, the Betti numbers of the ideal are computed and seen to be equivalent to Stanley's formula in the special case of face ideals of independence complexes of matroids.
 \end{abstract}

 \maketitle
 \section{Introduction}
   One of the foundations of oriented matroid theory is the topological representation theorem of Folkman and Lawrence \cite{FL}.   It says that  an oriented (simple) matroid can be realized uniquely as an arrangement of pseudospheres.  That there is no similar interpretation for the class of all matroids has been taken for granted.  For instance,``A non-coordinatizable matroid of abstract origin may be thought of as a geometric object only in a purely formal way, whereas an oriented matroid may always be thought of as a geometric-topological configuration on the $d$-sphere (or in projective space)'' \cite[pg.\ 19]{BLSWZ}.  Our main theorem is that the class of geometric lattices, which is cryptomorphic to the category of simple matroids,  is the same as the class of intersection lattices of arrangements of homotopy spheres.  

The interpretation of a geometric lattice as an arrangement of homotopy spheres is a natural generalization of the Folkman-Lawrence theorem. An oriented matroid realizable over $\R$ has a representation with geodesic spheres.  Allowing pseudospheres i.e., those which are homeomorphic,  but possibly not isometric to the unit sphere, leads to the category of (simple) oriented matroids.  If we further relax the conditions on the spheres to only homotopy equivalence to the standard sphere, then we are led to the  category of all (simple) matroids.

Some of the theory of oriented matroids which only depends on the underlying matroid can be extended to homotopy sphere arrangements.  Zaslavsky's enumerative theory for pseudosphere arrangements can be extended to the homotopy setting.  As in the oriented matroid representation theorem, the arrangement can be forced to be antipodal, so  a realization as homotopy {\it projective} spaces is also possible.  The minimal cellular resolutions of orientable matroid ideals developed in \cite{NPS} can be extended to arbitrary  matroids by using arrangements of homotopy spheres.

Our point of view is primarily through the lens of oriented matroids.  Hence the homotopy spheres which represent the atoms of the geometric lattice have codimension one.
In the future we hope to examine the point of view of  complex hyperplane arrangements and consider spheres of even codimension.

  The  matroid theory we require is in section \ref{matroids}, while section \ref{Steiner} presents  matroid Steiner complexes. We review a very general theory of arrangements of topological subspaces due to Ziegler and $\check{Z}$ivaljevi\'{c} in section \ref{subspaces}.  Arrangements of homotopy spheres and the representation theorems are in sections \ref{HSA} and \ref{main} respectively.  The last section extends the work of \cite{NPS} on minimal resolutions of face ideals of  independence complexes of matroids to matroid Steiner ideals.

Unlike some authors our homotopy spheres will not in general be manifolds.  The rest of our notation from topology is standard.  The join of two spaces $X$ and $Y$ is $X \ast Y,$ their one-point union, or wedge sum is $X \vee Y,$ while $X \simeq Y$ denotes homotopy equivalence.   
There are three facts from topology that we will use repeatedly without specific reference.  They are immediate consequences of well known theorems of Hurewicz and Whitehead.  See, for instance, \cite{Sp} for proofs. Let $\Gamma$ be CW-complex.  

\begin{itemize}
  \item
    If $\Gamma$ is acyclic and simply connected, then $\Gamma$ is contractible.
    
  \item
    If $\Gamma$ is simply connected and $\tilde{H}_0 (\Gamma) = \dots = H_{i-1}(\Gamma) = \{ 0 \},$ then $H_i(\Gamma) \cong \pi_i(\Gamma).$

  \item
    If $\Gamma$ and $\Gamma^\prime$ are homotopy spheres and $f:\Gamma \to \Gamma^\prime$ is a continuous map such that $f_\star:H_\star(\Gamma) \to H_\star(\Gamma^\prime)$ is an isomorphism, then $f$ is a homotopy equivalence.

\end{itemize}

\section{Matroids} \label{matroids}

In this section we give the basic definitions and results from matroid
theory that we will require.    Matroid definitions and notation are
as in \cite{O}.  Geometric semilattices are covered in \cite{WW}.  The characteristic polynomial and
M\"{o}bius invariant  can be found in \cite{Za}.  The beta invariant was introduced by Crapo \cite{C2}.

 There are numerous cryptomorphic definitions of matroids.  For us a {\it matroid} $M$ is a pair $(E, \mathcal{I}),\ E$ a non-empty
  finite set and $\mathcal{I}$ a distinguished set of subsets
  of $E.$ The members of $\mathcal{I}$ are called the {\it
  independent} subsets of $M$ and are required to satisfy:

\begin{itemize}
  \item The empty set is in $\mathcal{I}.$ 
  
  \item If $B$ is an independent set and $A \subseteq B,$ then $A$ is
    an independent set.  
   
  \item If $A$ and $B$ are independent sets such that $|A| < |B|,$
    then there exists an element $x \in B-A$ such that $A \cup \{x\}$ is
    independent.
\end{itemize} 
Matroid theory was introduced by Whitney \cite{W}.   The prototypical
example of a matroid is a  finite subset of a vector space over a field $\K$
with the canonical independent  sets. 
  Another
source of matroids is graph theory.  The cycle matroid of a graph is
the matroid whose finite set is the edge set of the graph and whose
independent sets are the acyclic subsets of edges. Most matroid
terminology can be traced back to these two types of examples.   

 The {\it circuits} of a matroid are its minimal
dependent sets.  If every circuit has cardinality at least three, then the matroid is {\it simple}.  A maximal
independent set is called a {\it basis}, and a subset which contains a basis is a {\it spanning} subset.  The maximal non-spanning subsets are the {\it hyperplanes} of $M.$  Every
basis of $M$ has the same cardinality.  The {\it rank} of $M,$ or
$r(M),$ is that common cardinality.  The {\it deletion} of $M$ at $e$ is 
denoted $M - e.$ It is the matroid whose
ground set is $E - \{e$\} and whose independent sets are simply those
members of $\mathcal{I}$ which do not contain $e.$  The {\it
contraction} of $M$ at $e$ is denoted $M/e.$  It is a matroid whose
ground set is also $E - \{e\}.$  If $e$ is a dependent element of $M$ then
$M/e = M - e.$  Otherwise, a subset $I$ of $E-\{e\}$ is independent in
$M/e$ if and only if $I \cup \{e\}$ is independent in $M.$  Deletion and
contraction for a subset $A$ of $E$ is defined by repeatedly deleting
or contracting each element of $A.$ The restriction of $M$ to $A$ is
$M-(E-A)$  which we shorten to $M|A$ or frequently just $A.$ The {\it rank} of
$A$ is $r(M|A)$ which we shorten to $r(A).$ Note that $r(\emptyset)=0.$ 

The {\it dual} of $M$ is $M^\star.$ It is the matroid whose ground set
is the same as $M$ and whose bases are the complements of the bases of
$M.$  For 
example, $U_{r,n}$ is the matroid defined by $E=\{1,2,\dots,n\}$ and 
$\mathcal{I}=\{A \subseteq E: |A| \le r\}.$  So, $U_{r,n}^\star = U_{n-r,n}.$  A circuit of $M^\star$ is a {\it cocircuit} of $M.$  The complement of a cocircuit is a hyperplane.  

The {\it free extension} of $M$ is $F(M).$  It is the matroid with ground 
set  $\tilde{E}(M) = E(M) \cup \{\tilde{e}\},$ where $\tilde{e} \notin E(M),$ and whose independent subsets are all  $A \subseteq \tilde{E}(M)$ such that $|A| \le r(M)$ and $A - \{\tilde{e}\}$ is independent in $M.$  The {\it free coextension} of $M$ is $(F(M^\star))^\star.$

A subset of $M$ is {\it closed} if adding any element to the subset
increases its rank.  The closed subsets of $M$ are also called {\it flats.} Examples of flats include the hyperplanes, also called the {\it coatoms} of $M.$ The closed subsets of $M$ with their inherited rank function form a ranked partially
ordered set under inclusion which we denote by  $L(M).$  Given two flats $X$ and $Y$ in $L(M)$ their {\it meet} is $X \wedge Y = X \cap Y.$ The meet of two flats is also their greatest lower bound in $L(M).$  The least upper bound of $X$ and $Y$ is their {\it join}, $X \vee Y,$ and is equal to $X \cup Y \cup \{e: \mbox{ there is a circuit } C,\  e \in C, C-\{e\} \subseteq X \cup Y\}.$ When the elements of a poset are topological spaces we rely on context to clarify whether $X \vee Y$ is their poset join or one-point union.  

\begin{defn}
  A finite ranked poset $L$  is a {\bf geometric lattice} if
   \begin{itemize}
     \item[(a)]
     $L$ is a lattice i.e., every pair of elements have a greatest lower bound and least upper bound.  In particular $L$ has a least element $\hat{0}$ and greatest element $\hat{1}.$

     \item[(b)]
     Every element of $L$ other than $\hat{0}$ is the join of atoms of $L.$

     \item[(c)]
     The rank function $r$  is semimodular, $r(X) + r(Y) \ge r(X \wedge Y) + r(X \vee Y).$

   \end{itemize}
\end{defn}

\begin{prop} \cite{O} \label{gl}
  If $M$ is matroid, then $L(M)$ is a geometric lattice.  Conversely, suppose $L$ is a geometric lattice.  Let $E$ be the atoms of $L.$  For any $A \subseteq E$ let $\bigvee A$ be the (poset) join of all the atoms in $A.$  Then $\mathcal{I} = \{A \subseteq E: r(\bigvee A) = |A|\}$ are the independent subsets of a matroid $M$ such that $L=L(M).$  
\end{prop}

The lower interval $[\hat{0}, X] \subseteq L(M)$ is the geometric lattice of $M|X.$  An upper interval $[X,\hat{1}]$ is isomorphic to the the geometric lattice of $M/X.$  For this reason $[X,\hat{1}]$ is frequently denoted $L/X.$ 

A {\it pointed geometric lattice} is a pair $(L,e)$ with $e$ a specified atom of $L.$ An isomorphism of pointed geometric lattices $\phi:(L,e) \to (L^\prime,e^\prime)$ is a lattice isomorphism  such that $\phi(e) = e^\prime.$ When $L = L(M)$ we will write $L(M,e)$ for the {\it geometric semilattice} $\{ X \in L: e \nleq X\}.$ (Geometric semilattices are called generalized affine matroids in \cite{McN}.) Even with the addition of a maximum element $L(M,e)$ is usually not a geometric lattice.

 Given any locally finite poset $L$ the  {\it M\"{o}bius function} on $L$
is the function $\mu:L \times L \to \Z$ which satisfies:
\begin{enumerate}
\item[(a)]
  If $X$ and $Y$ are incomparable elements of $L,$ then $\mu(X,Y) = 0.$

\item[(b)]
  For any $X$ in $L,\ \mu(X,X)=1.$

\item[(c)]
  For any $X,Y$ in $L, X < Y, \displaystyle{\sum_{X \le Z \le Y}
  \mu(X,Z) = 0.}$ 
\end{enumerate}

The {\it characteristic polynomial} of a geometric lattice $L$ is

$$p(L;t) = \sum_{X \in L} \mu(\hat{0},X) \ t^{r(\hat{1}) - r(X)}.$$

The {\it beta} invariant of a geometric lattice is
$$\beta(L) = (-1)^r \sum_{X \in L} r(X) \mu(\hat{0},X).$$

 \section{Matroid Steiner complexes} \label{Steiner}
 Matroid Steiner complexes were introduced in \cite{CP}.  We follow the presentation in \cite{Ch3}.
Let $M$ be a  matroid and $e \in E.$  The {\it port} of $M$ at $e$ is the set
$$ \mathcal{P} = \{C - \{e\}: C \mbox{ is a circuit of } M\}.$$

\noindent The matroid Steiner complex of $(M,e)$ is 

$$ \S(M,e) = \{F \subseteq E-\{e\}: \forall P \in \mathcal{P},  P \nsubseteq F \}.$$

Independence complexes of matroids are a subclass of matroid Steiner complexes. If $M$ is a matroid, then the independence complex of $M$ is 

$$\Delta(M) = \{F \subseteq E: F \mbox{ is independent} \}.$$

\noindent  Given a matroid $M$ let $\tilde{M}$ be the free coextension of 
$M$ with $\tilde{e}$ the extra point.  Then $\S(\tilde{M},\tilde{e}) = \Delta(M).$

\begin{thm} \cite{Ch3} \label{steiner}
  Let $e \in E, \ M$ a rank $r$  matroid.  Then
  
  $$\S(M,e) \simeq \bigvee^{\beta(M)}_{i=1} S^{r-2}.$$
\end{thm} 

\noindent Note that if $\beta(M)$ is zero, then $\S(M,e)$ is contractible.

 \section{Arrangements of subspaces} \label{subspaces}
 
 The theory of arrangements of subspaces of a topological space as presented in \cite{ZZ} will play a large part in our theory.  Here we present a CW-version of this theory.

Let $\A = \{A_1,\dots,A_m\}$ be a set of distinct subcomplexes of a finite CW-complex $\Gamma.$  Assume that $\A$ is closed under intersections.  Let $P$ be the poset $(\A,\le)$ where $A_i \le A_j$ if and only if $A_j \subseteq A_i.$  For any $p \in P, \Delta(P_{<p})$ is the order complex of the poset $P_{<p}=\{q \in P: q<p\}.$

\begin{thm} \cite{ZZ} \label{zz}
  Let $\V = \bigcup^m_{i=1} A_i.$  Suppose that for every $A_j < A_i,$ the inclusion map $A_i \hookrightarrow A_j$ is null-homotopic. Also assume that if $A_j$ has multiple components, then all non-empty images $A_i \hookrightarrow A_j$ are contained in the same component. Then

\begin{equation} \label{zz2}
 \V \simeq \bigvee^m_{i=1} (\Delta(P_{p<A_i}) \ast A_i) .
\end{equation}

\end{thm}
\noindent Note that $X \ast \emptyset = X.$

 \section{Arrangements of homotopy spheres} \label{HSA}
 
 A {\it homotopy $d$-sphere} is a $d$-dimensional CW-complex which is homotopy equivalent to  $S^d.$  It is convenient to let the empty set be a homotopy $(-1)$-sphere.

\begin{defn}
  A {\bf $d$-arrangement of homotopy spheres } consists of a $d$-dimensional homotopy sphere $S$ and a finite set of subcomplexes $ \A= \{S_1, \dots, S_n\}$ of $S$ each of which is a homotopy $(d-1)$-sphere.  In addition, 

  \begin{enumerate}
    \item[(a)]
       Every intersection of homotopy spheres in $\mathcal{A}$ is a homotopy sphere.

    \item[(b)]
      If $X$ is an intersection in $\A$ which is a $d^\prime$-dimensional homotopy sphere and $X \nsubseteq S_j,$ then $X \cap S_j$ is a $(d^\prime-1)$-dimensional homotopy sphere.

  \end{enumerate}
  
\end{defn}

Arrangements of homotopy spheres are a natural generalization of the arrangements of pseudospheres associated with oriented matroids. Many notions from pseudosphere arrangements can be carried over to the homotopy sphere case.  The {\it link} of $\A$ is $\V = \cup^n_{j=1} S_j.$ If $X$ is a non-empty  intersection in $\A,$ the {\it contraction} of $\A$ to $X$ is $\A/X.$ It is the $(\dim X)$-arrangement of homotopy spheres defined by letting $X$ be the ambient homotopy sphere and setting $\A/X$ equal to the collection of intersections $X \cap S_j, X \nsubseteq S_j.$  The {\it deletion} $\A - S_j$ is the arrangement $\{S_1,\dots,\hat{S}_j,\dots,S_n\}.$  We call $\A$  {\it essential} if $\cap^n_{j=1} S_j = \emptyset.$ As in the oriented matroid case, the intersection lattice of $\A$ plays a key role.

\begin{defn}
  Let $\A$ be an arrangement of homotopy spheres.  The {\bf intersection lattice} of $\A$ is the poset $L(\A)$ of intersections of elements of $\A$ ordered by reverse inclusion.  By convention, $S \in L(\A)$ as the empty intersection.
\end{defn}

As usual, $L(\A)/S_j \cong L(\A/S_j).$
We omit the elementary proof of the following.

\begin{prop}  Let $\A$ be a $d$-arrangement of homotopy spheres. Then $L(\A)$ is a  geometric lattice with rank function $r(X) = d- \dim X.$   If $\A$ is essential, then the rank of $L(\A)$ is $d+1.$
\end{prop}

In section \ref{main} we will prove the converse:  every rank $d+1$ geometric lattice is isomorphic to $L(\A)$ for some essential $d$-arrangement of homotopy spheres.  As with their pseudosphere counterparts, the homotopy type of the link of $\A$ only depends on $L(\A).$

\begin{prop} 
  If $\A$ is a $d$-arrangement of homotopy spheres, then
  $$\V \simeq \bigvee^{|p(L(\A);-1)|-1}_{i=1} S^{d-1}.$$
\end{prop}

\begin{proof}  We apply \ref{zz} to all of the subcomplexes in $L(\A)$ other than $S.$  So let $X \in L(\A) , X \neq S$  and let $r(X)$ be the rank of $X$ in $L(\A).$ Then $\Delta(L(\A)_{Y<X})$ is the order complex of the geometric lattice $[S,X] \subseteq L(\A).$ Hence \cite{Fo},

$$\Delta(L(\A)_{Y<X}) \simeq \bigvee^{|\mu(S,X)|}_{i=1} S^{r(X)-2}.$$

Since $X \simeq S^{d-r(X)}$ Theorem \ref{zz} implies

$$\V \simeq \bigvee_{X > S} \bigvee^{|\mu(S,X)|}_{i=1} S^{d-r(X)} \ast S^{r (X)-2}  = \bigvee^{|p(L(\A);-1)|-1}_{i=1} S^{d-1}.$$
 \end{proof}

In contrast to pseudosphere arrangements, an arrangement of homotopy spheres may not have a ``natural'' cell structure.  However, if $S$ does have a CW-structure which is related to $\A,$ then many of the enumerative invariants which describe the cell decomposition of a pseudosphere arrangement still hold in this more general situation.

\begin{figure}
  \begin{picture}(360,200)(-200,-100)
  \put (-190,0){\circle*{3}}
  \put (-150,0){\circle*{3}}
  \put (-90,0){\circle*{3}}
  \put (-50,0){\circle*{3}}
  \put (20,0){\circle*{3}}
  \put (60,0){\circle*{3}}
  \put (120,0){\circle*{3}}
  \put (160,0){\circle*{3}}
  \put (-120,50){\circle*{3}}
  \put (-120,-50){\circle*{3}}
  \put (90,50){\circle*{3}}
  \put (90,-50){\circle*{3}}
  \put (140,20){\circle*{3}}
  \put (-190,0){\line(1,0){40}}
  \put (-150,0){\line(3,5){30}}
  \put (-90,0){\line(-3,5){30}}
  \put (-90,0){\line(1,0){40}}
  \put (60,0){\line(3,5){30}}
  \put (120,0){\line(-3,5){30}}
  \put (20,0){\line(1,0){40}}
  \put (-150,0){\line(3,-5){30}}
  \put (-90,0){\line(-3,-5){30}}
  \put (120,0){\line(1,0){40}}
  \put (60,0){\line(3,-5){30}}
  \put (120,0){\line(-3,-5){30}}
  \put (160,0){\line(-1,1){20}}
  \put (90,50){\line(1,0){30}}
  \put (120,50){\circle*{3}}
  \put (-190,-20){$x_1$}
  \put(-50,-20){$y_1$}
  \put(-160,-20){$x_2$}
  \put(-83,-20){$y_2$}
  \put(-125,60){$x_3$}
  \put(-127,-61){$y_3$}
  \put(-170,40){$\A_1$}
  \put (20,-20){$x_1$}
  \put(50,-20){$x_3$}
  \put(125,-18){$x_2$}
  \put(155,-18){$y_2$}
  \put(88,60){$y_1$}
  \put(88,-61){$y_3$}
  \put(40,30){$\A_2$}
\end{picture}  
  \caption{Two representations of the three-point line where  $S_j=\{x_j,y_j\}.$} \label{3 point}
\end{figure}
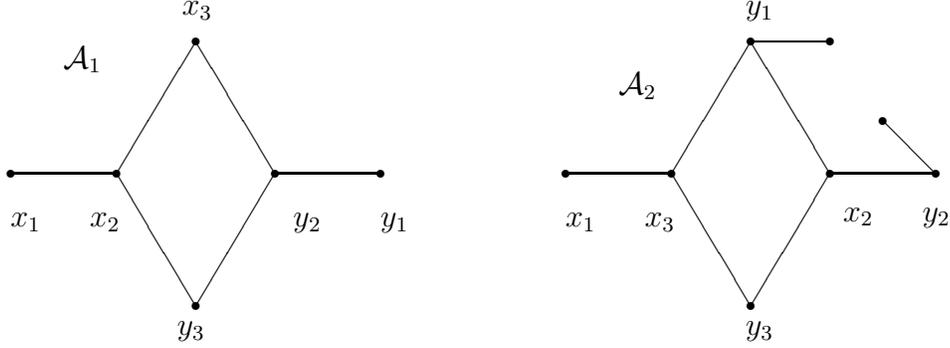  

\begin{example} \label{ex1}
  In figure \ref{3 point} both $\A_1$ and $\A_2$ are essential 1-arrangements of homotopy spheres with intersection lattice isomorphic to the three-point line. Only $\A_1$ has the same number of cells in each dimension as an essential pseudosphere arrangement with the same intersection lattice.  
\end{example}

\begin{defn}
  An arrangement of homotopy spheres is {\bf partitioned} if the $(d-1)$-skeleton of $S$ is contained in $\V.$  If every contraction of $\A$ is partitioned, then $\A$ is {\bf fully partitioned.}
\end{defn}

The CW-structure induced by a pseudosphere arrangement is  always fully partitioned.
Neither the contraction nor the deletion of a partitioned arrangement need be partitioned.  For instance, $\A_2$ in Figure \ref{3 point} is the deletion of a partitioned homotopy sphere arrangement with the intersection lattice of the 4-point line.

\begin{prop} \label{cells}
  If $\A$ is a partitioned  $d$-arrangement of homotopy spheres, then $S$ has $|p(L(\A);-1)| \ d$-dimensional cells.  
\end{prop}

\begin{proof}
 The proof is virtually identical to  the corresponding statements in \cite{Za2}.  For each cell $c$ in $S$ let $\psi(c) = \max \{X \in L(\A): c \subseteq X\}.$  For $X \in L(\A)$ define 

$$\tau(X) = \displaystyle\sum_{c \in \psi^{-1}(X)} (-1)^{\dim c}.$$

\noindent Note that we consider the empty set to be a cell of dimension minus one. Then, since each $X$ is a homotopy $(d-r(X))$-sphere, 

$$\sum_{Y \ge X} \tau(Y) = \chi(X) - 1 = (-1)^{d-r(X)}.$$

\noindent  M\"{o}bius inversion implies that,

$$\tau(X) = \sum_{Y \ge X} \mu(X,Y) (-1)^{d-r(X)}.$$

\noindent  Since $\A$ is partitioned, $\tau(S)$ is $(-1)^d$ times the number of $d$-dimensional cells.  
\end{proof}

\begin{cor}
  If every contraction $\A/X, r(X) \le i,$ is partitioned, then the number of $(d-i)$-dimensional cells is

$$\sum_{\stackrel{r(X) = i}{X \le Y}} |\mu(X,Y)|.$$

\end{cor}

\begin{proof}
  Under these conditions each $(d-i)$-cell is in exactly one rank-$i$ flat and $\tau(X)$ is $(-1)^{d-r(X)}$ times the number of $(d-i)$-cells in $X.$
\end{proof}

Unlike pseudosphere arrangements, $S-S_j$ need not consist of two contractible components.  Yet, under certain conditions, it is possible to recover Zaslavsky's enumerative results for the complex of bounded cells in a pseudosphere arrangement.

\begin{defn}
  An essential arrangement of homotopy spheres $\A$ is {\bf regular} with respect to $S_j \in \A$ if:
  
  \begin{itemize}
      \item
         For every $X \in L(\A)$ with $X \nsubseteq S_j,$ the subcomplex of $X$ generated by the vertices in $X - S_j$ consists of two contractible components. 

       \item
         For every coatom $Y > X$ such that $X \nsubseteq S_j$ the two vertices of $Y$ are in different components of $X - S_j.$
   \end{itemize}
\end{defn}

When $\A$ is regular with respect to $S_j$ a {\it bounded subcomplex} of $(\A,S_j)$ is one of the two components of the subcomplex generated by the vertices in $S - S_j.$  In this situation we define,

$$L(\A,S_j) = \{X \in L(\A): X \nsubseteq S_j\}.$$

\begin{prop}
  Suppose $\A$ is a partitioned $d$-arrangement of homotopy spheres which is regular with respect to $S_j.$  Then the number of $d$-cells of a bounded subcomplex of $(\A,S_j)$ is

\begin{equation} \label{bounded cells}
|\displaystyle\sum_{X \in L(\A,S_j)} \mu(S,X)| = \beta(L(\A)).
\end{equation}

\end{prop}

\begin{proof}  
    Let $B$ be a bounded subcomplex of $(\A,S_j).$  For each $X$ in $L(\A,S_j)$ let $X_B = B \cap X.$  For each cell $c \subseteq B,$ define $\psi(c)$ to be the maximal $X \in L(\A,S_j)$ such that $c \subseteq X_B.$  Define 

$$\tau(X) = \displaystyle\sum_{\stackrel{c \in \psi^{-1}(X)}{c \neq \emptyset}} (-1)^{\dim c}.$$   

\noindent Since each $X_B$ is contractible its Euler characteristic is one.  Therefore, for any $X \in  L(\A,S_j)$

$$ \sum_{\stackrel{Y \ge X}{Y \in L(\A,S_j)}} \tau(Y) = 1.$$

\noindent  M\"{o}bius inversion and  the fact that $\tau(S)$ is $(-1)^d$ times the number of $d$-cells in $B$ imply that the number of $d$-cells is the left-hand side of (\ref{bounded cells}).  That this equals $\beta(L(\A))$ is  \cite[pg.\ 77]{Za2}.  
\end{proof}

\begin{cor} \label{bounded f-vector}
  If $\A$ is a fully partitioned $d$-arrangement of homotopy spheres which is regular with respect to $S_j,$ then the number of $(d-i)$-cells in a bounded subcomplex of $S$ is

$$|\sum_{\stackrel{Y \in L(\A,S_j)}{r(X) = i }} \mu(X,Y)|$$
\end{cor}

As is evident from the proofs, the enumerative results in this section only depend on the Euler characteristic of spheres and contractible spaces.  Hence these results would apply to arrangements of any spaces with the same Euler characteristics.  This idea is explored in much greater generality in \cite{Za3}.

Not every invariant which only depends on the underlying matroid carries over from  pseudosphere arrangements to homotopy sphere arrangements.  The flag $f$-vector of a pseudosphere arrangement only depends on the intersection lattice \cite{BSt2}.  However, as the following example shows, this is not true for homotopy sphere arrangements even if the arrangement is fully partitioned and a combinatorial CW-complex.

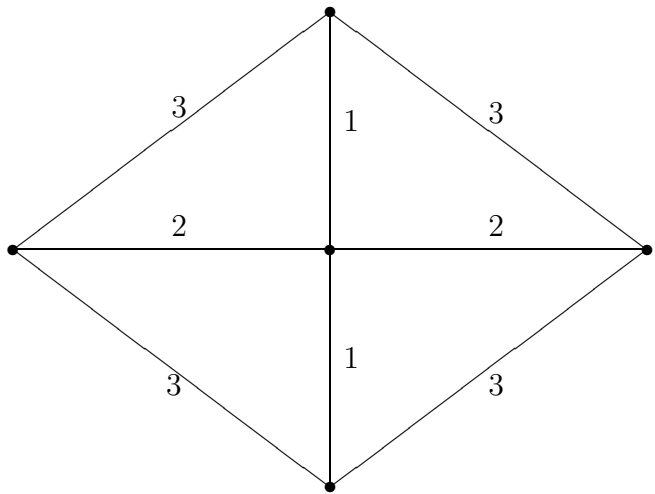
\begin{figure}
  \begin{picture}(300,200)(-150,-100)  
  \put(-120,0){\circle*{4}}
  \put(0,0){\circle*{4}}
  \put(120,0){\circle*{4}}
  \put(0,90){\circle*{4}}
  \put(0,-90){\circle*{4}}

  \put(-120,0){\line(4,3){120}}
  \put(0,90){\line(4,-3){120}}
  \put(-120,0){\line(4,-3){120}}
  \put(0,-90){\line(4,3){120}}
  \put(-120,0){\line(1,0){240}}
  \put(0,90){\line(0,-1){180}}
  
  \put(-60,5){2}
  \put(60,5){2}
  \put(5,45){1}
  \put(5,-45){1}
  \put(-60,50){3}
  \put(60,48){3}
  \put(60,-55){3}
  \put(-62,-55){3}
\end{picture}  
  \caption{1-skeleton of the upper hemisphere of the Boolean algebra with three atoms. } \label{flag-vector}
\end{figure}

\begin{example}
  Figure \ref{flag-vector} shows the 1-skeleton of the upper hemisphere for a homotopy 2-arrangement of the Boolean algebra with three atoms.  If we attach four interior triangles and glue a copy of this 2-complex along the outside boundary, then the resulting arrangement has twenty-four $0$-cell $\subset 2$-cell incidences.  Alternatively, we could attach one square along the outside boundary and three interior triangles and then glue a copy of this complex along the outside boundary.  Now there are twenty-six $0$-cell $\subset 2$-cell incidences.  
\end{example}

\section{The representation theorem} \label{main}

\begin{thm} \label{rep thm}
  Let $(L,e)$ be a rank $r$ pointed geometric lattice.  There exists a fully partitioned  essential $(r-1)$-arrangement of homotopy spheres $\A$ which is regular with respect to $S_1$ such that $(L,e)$ and $(L(\A),S_1)$ are isomorphic as pointed lattices.  Furthermore, the arrangement can be constructed so that there exists a fixed-point free involution of $S$ which preserves $\A.$
\end{thm}

\begin{proof} 
  The proof is by induction on $r.$   It will be evident from the construction that $(\A,S_1)$ has the following additional property.  Let $a_1,\dots,a_m$ be the coatoms of $L,$ and let $\{x_1,y_1\},\dots,\{x_m,y_m\}$ be the corresponding zero-spheres in $\A.$   Then the map which takes $X \in L$ to the subcomplex induced by the vertices $\{x_i,y_i:X \le a_i\}$ is a lattice isomorphism.

We begin with $r=2.$ While simply putting pairs of antipodal points around a circle satisfies the theorem, we prefer to give a procedure which produces all possible arrangements which satisfy the theorem as it is indicative of what happens in higher ranks.  Since $L$ is rank two it consists of $\hat{0}, \hat{1}$ and coatoms (=atoms) $\{a_1,\dots,a_m \},$ where $a_1 = e.$   So $S$ has $2m$ vertices which we label $x_1,y_1,\dots,x_m,y_m$  where the pair $\{x_i,y_i\}$ corresponds to the atom $a_i.$  Choose any spanning tree $D_b$ of $\{x_2,\dots,x_m\}$ and extend it to a spanning tree $D$ of $\{x_1,\dots,x_m,y_1\}.$  Let $D^\prime$ be the mirror image of $D$ with $\{y_1,\dots,y_m,x_1\}$ replacing $\{x_1,\dots,x_m,y_1\}.$   Finally, let $S$ be $D \cup D^\prime.$  Now, $S$ deformation retracts to the circle formed by the unique path from $x_1$ to $y_1$ in $D$ concatenated with the unique path from $y_1$ to $x_1$ in $D^\prime.$  In addition, $\A = \{\{x_1,y_1\},\dots,\{x_m,y_m\}\}$ is a 1-arrangement of homotopy zero-spheres which satisfies the theorem, where the involution is the  one induced by switching $x_i$ and $y_i.$  We also note that any arrangement which satisfies the theorem must be of this form. Figure \ref{rank 2} shows all three stages of this construction for the five-point line.

\begin{figure} 
  \begin{picture}(300,200)(-200,-100)  
  \put (-100,100){$x_1$}
  \put (-100,-100){$y_1$}
  \put (-190,59){$x_2$}
  \put (-130,0){$x_3$}
  \put (-120,56){$x_5$}
  \put (-192,-55){$x_4$}
  \put (-97,95){\circle*{4}}
  \put (-97,-89){\circle*{4}}
  \put (-175,-52){\circle*{4}}
  \put (-175,58){\circle*{4}}
  \put (-127,57){\circle*{4}}
  \put (-127,10){\circle*{4}}
  \put (-175,58){\line(0,-1){110}}
  \put (-175,58){\line(1,-1){48}}
  \put (-127,10){\line(0,1){47}}
  \put (-140,-55){\Large{$D_b$}}
  
  \put (50,100){$x_1$}
  \put (50,-100){$y_1$}
  \put (-40,59){$x_2$}
  \put (20,0){$x_3$}
  \put (30,56){$x_5$}
  \put (-42,-55){$x_4$}
  \put (53,95){\circle*{4}}
  \put (53,-89){\circle*{4}}
  \put (-25,-52){\circle*{4}}
  \put (-25,58){\circle*{4}}
  \put (23,57){\circle*{4}}
  \put (23,10){\circle*{4}}
  \put (-25,58){\line(0,-1){110}}
  \put (-25,58){\line(1,-1){48}}
  \put (23,10){\line(0,1){47}}
  \put (53,95){\line(-2,-1){80}}
  \put (-25,-52){\line(2,-1){78}}
  \put (10,-55){\Large{$D$}}
\end{picture}  
  \begin{picture}(300,200)(-200,-100)

  \put (-50,100){$x_1$}
  \put (-50,-100){$y_1$}
  \put (-140,59){$x_2$}
  \put (-80,0){$x_3$}
  \put (-70,56){$x_5$}
  \put (-142,-55){$x_4$}
  \put (-47,95){\circle*{4}}
  \put (-47,-89){\circle*{4}}
  \put (-125,-52){\circle*{4}}
  \put (-125,58){\circle*{4}}
  \put (-77,57){\circle*{4}}
  \put (-77,10){\circle*{4}}
  \put (-125,58){\line(0,-1){110}}
  \put (-125,58){\line(1,-1){48}}
  \put (-77,10){\line(0,1){47}}
  \put (-47,95){\line(-2,-1){80}}
  \put (-125,-52){\line(2,-1){78}}

  \put (25,63){$y_4$}
  \put (-30,20){$y_3$}
  \put (-30,-52){$y_5$}
  \put (35,-55){$y_2$}
  \put (25,-52){\circle*{4}}
  \put (25,58){\circle*{4}}
  \put (-23,-36){\circle*{4}}
  \put (-23,10){\circle*{4}}
  \put (25,58){\line(0,-1){110}}
  \put (25,-52){\line(-3,4){48}}
  \put (-23,10){\line(0,-1){47}}
  \put (25,-90){\Large{$S$}}
  \put (-47,95){\line(2,-1){73}}
  \put (25,-52){\line(-2,-1){73}}
  
\end{picture}  
  \caption{A 1-arrangement of homotopy spheres representing the 5-point line.} \label{rank 2}
\end{figure}
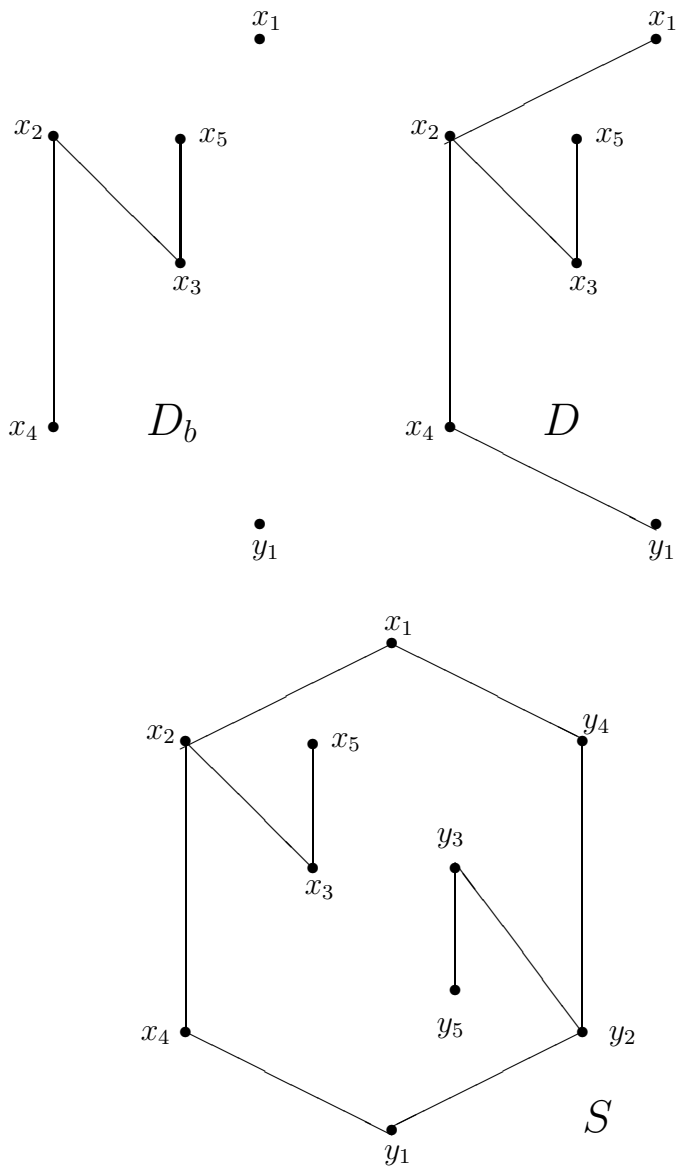

 Now suppose the rank of $L$ is three. Let $a_1,\dots,a_m$ be the coatoms of $L$ and let $e_1,\dots,e_n$ be the atoms of $L$ with $e_1=e.$ As above, we choose vertices $x_1,y_1,\dots,x_m,y_m$ with the pair $\{x_i,y_i\}$ corresponding to $a_i.$  The interval $[e_1,\hat{1}]$ in $L$ is a rank 2 geometric lattice. Therefore, we can use the construction above to obtain an arrangement of zero spheres on a homotopy 1-sphere $S_1$ which satisfies the theorem and represents $[e_1,\hat{1}]$ with vertices consisting of all the $x_i$ and $y_i$ such that $e_1 < a_i.$ 

Let $e_j$ be another atom of $L.$  Let $W_j = \{ x_i:e_j < a_i\} \cup \{y_i: a_i = e_1 \vee e_j\}.$  As above we choose a spanning tree $V_j$ on $W_j$ which is also a spanning tree on the vertices of $W_j$ which do not represent $e_1 \vee e_j.$  Let $\V$ be the union of $S_1$ and all the $V_j.$   Let $\V_b$ be the subcomplex of $\V$ generated by all vertices representing $\{a_i: e_1 \not< a_i.\}$  Since $\V_b$ is a connected graph, we can attach $2$-cells so that the resulting space, which we call $D_b,$ is contractible.  The CW-complex $\V \cup D_b$ is homotopy equivalent to a connected graph, so  we can attach 2-cells to it so that the resulting space is contractible.  Call this two-dimensional CW-complex $D.$  Note that none of the 2-cells attached to $\V \cup D_b$ have their boundary completely contained in $D_b.$  Otherwise $D$ would not be acyclic.  Let $D^\prime$ be the mirror image of $D$ obtained by switching the roles of $x_i$ and $y_i$ and let $S$ be the union of $D$ and $D^\prime$ glued along $S_1.$ The intersection $D \cap D^\prime = S_1,$ so $S$ is a homotopy $2$-sphere.  For each $i, 1 \le i \le n,$ let $S_i$ be the subcomplex of $S$ induced by the all vertices $x_j, y_j$ such that $e_i < a_j.$  By construction each $S_i$ is a homotopy one-sphere and the arrangement $\A = \{S_1, \dots, S_n\}$ satisfies the theorem.

\begin{example}
  Figure \ref{fano} shows one possible way of constructing $\V$ for $(L,1)$ where $L$ is the Fano plane as pictured.  The 1-cells which will be used to form $S_1,S_2$ and $S_3$ are also labeled. The subcomplex induced by $\{D,E,F,G\}$ is $\V_b$ and is the 1-skeleton of a tetrahedron.  Attaching triangles $EFG, DEF$ and $DFG$ is one of infinitely many ways to construct $\V \cup D_b.$
 \end{example}
 
 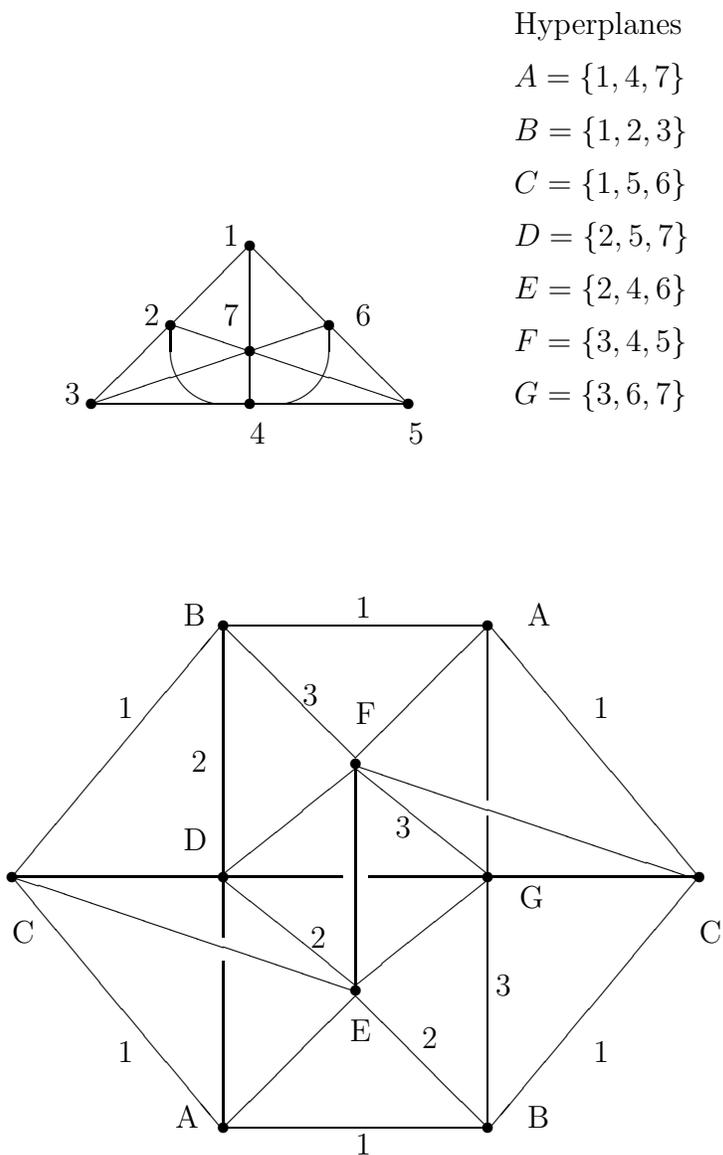
\begin{figure}
   \begin{picture}(300,200)(-150,-100)  
  \put(-50,0){\circle*{4}}
  \put(-80,-30){\circle*{4}}
  \put(-110,-60){\circle*{4}}
  \put(-50,-60){\circle*{4}}
  \put(10,-60){\circle*{4}}
  \put(-20,-30){\circle*{4}}
  \put(-50,-40){\circle*{4}}
  
  \put(-50,0){\line(-1,-1){60}}
  \put(-50,0){\line(1,-1){60}}
  \put(-110,-60){\line(1,0){120}}
  \put(-50,0){\line(0,-1){60}}
  \put(-110,-60){\line(3,1){90}}
  \put(10,-60){\line(-3,1){90}}
  \put(-50,-30){\oval(60,60)[bl]}
  \put(-50,-30){\oval(60,60)[br]}

  \put(-60,0){1}
  \put(-90,-30){2}
  \put(-120,-60){3}
  \put(-50,-75){4}
  \put(10,-75){5}
  \put(-10,-30){6}
  \put(-60,-30){7}
  
  \put (50,80){Hyperplanes}
  \put (50,60){$A = \{1,4,7\}$}
  \put (50,40){$B = \{1,2,3\}$}
  \put (50,20){$C = \{1,5,6\}$}
  \put (50,0){$D = \{2,5,7\}$}
  \put (50,-20){$E = \{2,4,6\}$}
  \put (50,-40){$F = \{3,4,5\}$}
  \put (50,-60){$G = \{3,6,7\}$}
\end{picture}  
   \begin{picture}(300,220)(-150,-110)  
  \put(-140,0){\circle*{4}}
  \put(-60,0){\circle*{4}}
  \put(40,0){\circle*{4}}
  \put(120,0){\circle*{4}}
  \put(-60,95){\circle*{4}}
  \put(-60,-95){\circle*{4}}
  \put(40,95){\circle*{4}}
  \put(40,-95){\circle*{4}}
  \put(-10,43){\circle*{4}}
  \put(-10,-43){\circle*{4}}
  \put(-140,0){\line(1,0){80}}
  \put(-60,0){\line(1,0){45}}
  \put(-5,0){\line(1,0){45}}
  \put(40,0){\line(1,0){80}}
  \put(-60,94){\line(0,-1){94}}
  \put(-60,0){\line(0,-1){23}}
  \put(-60,-32){\line(0,-1){63}}
  \put(-140,0){\line(3,-1){130}}
  \put(-10,42){\line(3,-1){130}}
  \put(40,95){\line(0,-1){66}}
  \put(40,23){\line(0,-1){23}}
  \put(-140,0){\line(5,6){80}}
  \put(-60,95){\line(1,-1){50}}
  \put(40,-95){\line(-1,1){50}}
  \put(-60,-95){\line(1,1){50}}
  \put(40,95){\line(-1,-1){50}}
  \put(-60,1){\line(5,4){50}}
  \put(-60,-1){\line(5,-4){50}}
  \put(40,-1){\line(-5,-4){50}}
  \put(40,1){\line(-5,4){50}}
  \put(-10,43){\line(0,-1){86}}
  \put(-140,0){\line(5,-6){80}}
  \put(120,0){\line(-5,6){80}}
  \put(120,0){\line(-5,-6){80}}
  \put(-60,95){\line(1,0){100}}
  \put(-60,-95){\line(1,0){100}}
  \put(40,0){\line(0,-1){95}}
  \put(-140,-25){C}
  \put(-75,10){D}
  \put(52,-12){G}
  \put(120,-25){C}
  \put(-75,95){B}
  \put(-78,-95){A}
  \put(55,95){A}
  \put(55,-95){B}
  \put(-10,58){F}
  \put(-12,-62){E}
  \put(-100,60){1}
  \put(-10,98){1}
  \put(-10,-105){1}
  \put(-100,-70){1}
  \put(80,-70){1}
  \put(80,60){1}
  \put(-72,40){2}
  \put(-27,-27){2}
  \put(15,-65){2}
  \put(43,-45){3}
  \put(5,15){3}
  \put(-30,65){3}

\end{picture}  
   \caption{The Fano plane} \label{fano}
 \end{figure}  
  
For the induction step, assume the rank of $L$ is greater than three.  As above, let $a_1,\dots,a_m$ be the coatoms of $L$ and let $e_1,\dots,e_n$ be the atoms of $L$ with $e_1 = e.$ In addition, let $x_1,y_1,\dots,x_m,y_m$ be vertices with the pair $\{x_i,y_i\}$ corresponding to $a_i.$  The interval $L_1 = [e_1,\hat{1}]$ in $L$ is a rank $r-1$ geometric lattice. Using the inductive algorithm we construct an arrangement $\A_1$ of homotopy $(r-2)$-spheres which satisfy the theorem for $(L_1,e^\prime),$ where $e^\prime$ is any atom of $L_1.$  

Now consider the codimension-2 flats of $L.$  Those which contain $e_1$ are already represented in $\A_1.$  For any other codimension-2 flat $X,$  we construct a spanning tree $V_X$ on $W_X=\{ x_i:X < a_i\} \cup \{y_i: a_i = X \vee e_1\}$ which is also a spanning tree on the vertices of $W_X$ which do not represent $X \vee e_1.$  

Proceeding inductively on the corank, for each corank-$k$ flat $X \neq \hat{0}, e_1 \nleq X$ we  construct a contractible space $V_X$ on the vertices $W_X=\{x_i: X < a_i\} \cup \{y_i:X \vee e_1 \le a_i\}$ by adding only $(k-1)$-dimensional cells to the complex determined by flats of lower corank.  In addition, we make sure that the subcomplex of $V_X$ generated by the vertices $\{x_i: e_1 \nleq a_i\}$ is also contractible.  Flats which lie above $e_1$ are already represented in $S_1$

Let $\V =\bigcup_{X \neq \hat{0}} V_X$ and let $\V_b$ be the subcomplex of $\V$ determined by the vertices whose corresponding coatoms do not lie above $e_1.$
Consider the covering of $\V_b$ by the subcomplexes $U_j =  V_{e_j} \cap \V_b, 2 \le j.$  By construction, all non-void intersections of the members of the cover are contractible so the classical nerve theorem applies.  Since a collection of the $U_j$ have void intersection if and only if the poset join of the corresponding atoms contains $e_1, \V_b$ is homotopy equivalent to 

$$
\Delta(L,e)  = \{F \subseteq \{e_2,\dots,e_m\}: e_1 \notin \bigvee F \} 
$$

 This complex is just  the matroid Steiner complex $\S(M,e)$ where $M$ is the matroid associated to $L$ as described in Proposition \ref{gl}. Hence, by Theorem \ref{steiner} it has the homotopy type of  a wedge of $\beta(L)$ spheres of dimension $(r-2).$    Now we attach $\beta(L)$ cells of dimension $r-1$ to $\V_b$ in any way which kills $H_{r-2}(\V_b)$ and call this space $D_b.$ Since $D_b$ is acyclic and simply connected it is contractible.  

Now consider the cover of  $\V$ with the subcomplexes $V_{e_j}$ and $S_1.$  Apply Theorem \ref{zz} to the arrangement of subcomplexes consisting of this cover and all of its intersections. Call this arrangement $\B.$  For $B \in \B$ define $\phi(B)$ to be the flat in $L$ which is the intersection of all of the hyperplanes whose corresponding vertices are incident to $B.$  Then $\phi$ is a lattice isomorphism from $\B$ to $L.$ The elements of $\B$ are of two types.  If $e_1 \notin \phi(B),$ then $B$ is contractible. Otherwise, if  $e_1 \le \phi(B), $ then $B$ is a homotopy $r-r(\phi(B))$-sphere.  As $\Delta(\B_{< B}) \simeq \Delta(L_{< \phi(B)}),$ every non-contractible term of (\ref{zz2}) is a wedge of $(r-2)$-spheres. 

Since $D_b$ is contractible and $D_b \hookrightarrow D_b \cup \V$ is a cofibration, the reduced homology of $D_b \cup \V$ is the same as the homology of the pair $(\V,\V_b).$  As $\V$ is an (r-2)-dimensional CW-complex $H_{r-2}(\V,\V_b)$ is torsion-free. This plus the long exact sequence of the pair implies that $(\V,\V_b)$, and hence $D_b \cup \V,$ has the homology of a wedge of $(r-2)$-spheres. Let $D$ be any $(r-1)$-dimensional CW-complex obtained by gluing $(r-1)$-dimensional cells to $D_b \cup \V$ so that the resulting space is acyclic.  As before, none of these cells have their boundary contained in $D_b.$  Since $\V$ is simply connected $D$ is simply connected, and hence contractible.  Let $D^\prime$ be the mirror image of $D$ induced by switching  $x_i$ with $y_i.$ Finally, let $S = D \cup D^\prime$ glued along $S_1.$ Since $S$ is the union of two contractible spaces whose intersection is a homotopy $(r-2)$-sphere, $S$ is a homotopy $(r-1)$-sphere.  Similar reasoning shows that for each flat $X \in L$ the subcomplex of $S$ generated by the vertices corresponding to the hyperplanes above $X$ in $L$ is a homotopy $(r-r(X))$-sphere and the intersection lattice of the arrangement of homotopy spheres determined by the atoms of $L$ is isomorphic to $L.$  The involution induced by switching $x_i$ and $y_i$ for each $i$ is fixed-point free and preserves the arrangement.  The construction insures that this arrangement is fully partitioned and is also regular with respect to $S_1.$
  
\end{proof}

%\begin{ques}
%  What conditions on $L$ allow the construction of $\A$ so that it is regular %with respect to all the atoms of $L?$
%\end{ques}

\begin{remark}
  The involution means that we can also develop a theory of homotopy projective space arrangements.
\end{remark}

\begin{defn}
  Let $\A = \{S_1,\dots,S_n\}$ and $\B  = \{T_1,\dots,T_n\}$ be $d$-arrangements of homotopy spheres.  Then $\A$ and $\B$ are {\bf homotopy equivalent arrangements} if there exists a lattice isomorphism $\phi: L(\A) \to L(\B)$ and homotopy equivalences $f:S \to T, g:T \to S$ such that:

\begin{enumerate}
  \item
     For each $X \in L(\A), f(X) \subseteq \phi(X),$ and  $f:X \to \phi(X)$ is a homotopy equivalence.
     
  \item
     For each $Y \in L(\B), g(Y) \subseteq \phi^{-1}(Y)$ and $g:Y \to \phi^{-1}(Y)$ is a homotopy equivalence.
     
\end{enumerate}
\end{defn}

\begin{thm}[Uniqueness] \label{uniqueness}
  Let $\A, \B$ be essential $d$-arrangements of homotopy spheres such that $L(\A) \cong L(\B).$  Then $\A$ and $\B$ are homotopy equivalent  arrangements.
\end{thm}

\begin{proof}
  Let $\phi: L(\A) \to L(\B)$ be a lattice isomorphism.  For every cell of $S$ define $\psi(c)$ to be $\max \{X \in L(\A): c \subseteq X\}.$ By symmetry it is sufficient to construct $f$ in the above definition.  We will build up $f$ by defining  maps $f_i$ on the $i$-skeletons of $S$ inductively which satisfy the following properties:

\begin{itemize}
  \item
      The restriction of $f_{i+1}$ to the $i$-skeleton is $f_i.$
      
  \item
      For each corank-$(i+1)$ flat $X \in L(\A), f_i:X \to  \phi(X)$ is a homotopy equivalence. 

  \item
      If $\dim c \le i,$ then $f_i(c) \subseteq \phi(\psi(c)).$
\end{itemize}
      
         Since the arrangements are essential the coatoms are all homeomorphic to two disjoint points. Choose a homeomorphism $h$ of the union of the coatoms of $\A$ to the coatoms of $\B$ which preserves $\phi.$  Now extend $h$ to  $f_0$ by arbitrarily choosing any image point in $\phi(\psi(v))$ for any vertex $v$ which is not in a coatom of $\A.$

Now assume that $f_{i-1}$ has been defined.
Let $X$ be a corank $(i+1)$-flat of $L(\A).$  Since $\A$ is essential, $X$ is a homotopy $i$-sphere.  Let $c$ be an $i$-cell of $X.$ The definition of $f_{i-1}$ insures that $f_{i-1}(\partial(c)) \subseteq \phi(X).$  As $\phi(X)$ is also a homotopy $i$-sphere, there is a map $f_c: \bar{c} \to \phi(X)$ such that $f_c$ equals $f_{i-1}$ when restricted to $\partial(c).$ Putting all of these maps together gives  a map $f_X: X \to \phi(X).$  The induced map in homology, $(f_X)_\star: H_i(X) \to H_i(\phi(X)),$ is multiplication by $n_X$ after choosing generators for the respective homology groups.  If $n_X = \pm 1 ,$ then $f_X$ is a homotopy equivalence.  If not, choose any $i$-cell $c$ in $X$ and redefine $f_c$ as follows.  Let $\alpha:(D^i,S^{i-1}) \to X$ be the attaching map for $c.$ Let $D^i(1/2)$  be the closed ball of radius one-half and let  $S^{i-1}(1/2) = \partial D^i(1/2).$ Replace $f_c$ with $\tilde{f}_c$ which satisfies: 

\begin{itemize}
  \item
  $\tilde{f}_c$ restricted to $S^{i-1}(1/2)$ is constant and equal to $f_c(0).$ The induced map $(\tilde{f}_c)_\star:H_i(D^i(1/2),S^{i-1}(1/2)) \to H_i(\phi(X))$ is multiplication by $1 - n_X$ with respect to the appropriate generators.
  \item 
   $\tilde{f}_c(\mathbf{x})= f_c((2|\mathbf{x}|-1) \mathbf{x})$ for $|\mathbf{x}| \ge 1/2.$
\end{itemize}

The new $f_X$ induces an isomorphism on homology and hence is a homotopy equivalence.  If an $i$-cell $c$  is not contained in any corank-$(i+1)$ flat, then $\psi(c)$ is at least $i$ connected, so we define an arbitrary map $f_c: c \to \phi(\psi(c))$ which equals $f_{i-1}$ on $\partial c.$   Putting all of the $f_c$ together gives the required map $f_i.$
\end{proof}

\section{Minimal cellular resolutions of matroid Steiner ideals}  \label{mcr}
  One approach to finding syzygies of monomial ideals is through minimal cellular resolutions.  The following presentation of minimal cellular resolutions is taken from \cite{NPS}.  Let $\K$ be a field and let $I$ be the monomial ideal $<m_1,\dots,m_s>$ in the polynomial ring $\mathbf{k}[x_1,\dots,x_n]$ which we denote by  $\mathbf{k[x]}.$   Let $\Gamma$ be a CW-complex with $s$ vertices $v_1,\dots,v_s$ which are labeled with the monomials $m_1,\dots,m_s.$  Write $c \ge c^\prime$ whenever a cell $c^\prime$ belongs to the closure of another cell $c.$  Label each cell $c$ of $\Gamma$ with the monomial $m_c = \lcm\{m_i:v_i \le c\},$ the least common multiple of the monomials labeling the vertices of $c.$  Set $m_\emptyset = 1$ for the empty cell of $\Gamma.$  The principal ideal $<m_c>$ is identified with the free $\N^n$-graded $\mathbf{k[x]}$-module of rank 1 with generator in degree $\deg m_c.$  For a pair of cells $c \ge c^\prime,$ let $p^{c^\prime}_c:<m_c> \hookrightarrow <m_{c^\prime}>$ be the inclusion map of ideals.  It is a degree-preserving homomorphism of $\N^n$-graded modules.

Fix an orientation of each cell in $\Gamma,$ and define the {\it cellular complex} $C_\star(\Gamma,I)$

$$ \dots \stackrel{\partial_3}{\longrightarrow} C_2 \stackrel{\partial_2}{\longrightarrow} C_1 \stackrel{\partial_1}{\longrightarrow} C_0 \stackrel{\partial_0}{\longrightarrow} C_{-1} = \mathbf{k[x]}$$

\noindent as follows. The $\N^n$-graded $\mathbf{k[x]}$-module of $i$-chains is

$$C_i = \bigoplus_{c: \dim c = i} <m_c>,$$

\noindent where the direct sum is over all $i$-dimensional cells $c$ of $\Gamma.$  The differential $\partial_i: C_i \to C_{i-1}$ is defined on the component $<m_c>$ as the weighted sum of the maps $p^{c^\prime}_c:$

$$\partial_i = \sum_{\stackrel{c^\prime \le c}{\dim c^\prime = i-1}} [c:c^\prime] \ p^{c^\prime}_c,$$

\noindent where $[c:c^\prime] \in \Z$ is the incidence coefficient of the oriented cells $c$ and $c^\prime$ in the usual topological sense.  The differential $\partial_i$ preserves the $\N^n$-grading of $\mathbf{k[x]}$-modules.  Note that if $m_1 = \dots = m_s =1,$ then $C_\star(\Gamma,I)$ is the usual chain complex of $\Gamma$ over $\mathbf{k[x]}.$  For any monomial $m \in \mathbf{k[x]},$ define $\Gamma_{\le m}$ to be the subcomplex of $\Gamma$ consisting of all cells $c$ whose label $m_c$ divides $m.$  We call any such $\Gamma_{\le m}$  an $I$-essential subcomplex of $\Gamma.$

\begin{prop} \cite[Proposition 1.2]{BSt}
  The cellular complex $C_\star(\Gamma,I)$ is exact if and only if every $I$-essential subcomplex of $\Gamma$ is acyclic over $\mathbf{k}.$  Moreover, if $m_c \neq m_{c^\prime}$ for any $c > c^\prime,$ then $C_\star(\Gamma,I)$ gives a minimal free resolution of $I.$
\end{prop}

If both of the conditions of the above proposition are met then we call $\Gamma$ an {\it $I$-complex} and $C_\star(\Gamma,I)$ a {\it minimal cellular resolution} of $I.$  Recall that $\beta_i(I)$ is the $\K$-dimension of the $i^{th}$ free module in a minimal free resolution of $I.$  When $\Gamma$ is an $I$-complex the number of $i$-dimensional cells in $\Gamma$ is $\beta_i(I).$  

Given an abstract simplicial complex $\Delta$ with vertices $v_1,\dots, v_n$ the {\it face ideal} of $\Delta$  in $\mathbf{k[x]}$ is

$$I_\Delta = <\{ x_{i_1} \cdots x_{i_s} : \{v_{i_1}, \dots, v_{i_s}\} \notin \Delta \}>.$$

When $\Delta$ is a matroid Steiner complex we call $I_\Delta$ a {\it matroid Steiner ideal.}
As pointed out it section \ref{steiner} independence complexes of matroids are a special subclass of matroid Steiner complexes.   The problem of finding minimal resolutions of $I_\Delta$ when $\Delta$ is the independence complex of a matroid $M$ was examined in \cite{NPS}.  When $M$ is an orientable matroid Novik et al. showed that the bounded subcomplex of any pseudosphere arrangement which realizes $M^\star$ extended by a free point is an $I_\Delta$-complex. 

Let $I$ be a matroid Steiner ideal.  The topological representation theorem allows a complete description of all possible equivalence classes of complexes which are $I$-complexes for every field $\K.$  Two $I$-complexes are {\it equivalent} if they have the same cellular resolution (up to orientation). Acyclic 2-complexes which are not simply connected show that it is possible for two equivalent $I$-complexes to be  homotopy inequivalent.

\begin{thm}  \label{I-complexes}
  Let $I = I_{\S(M,e)}$ be a matroid Steiner ideal.  Let $(\A,S_1),$ be a fully partitioned  arrangement of homotopy spheres which is regular with respect to $S_1$ such that $(L(\A), S_1) \cong (L(M^\star),e).$  Then a bounded subcomplex of $(\A,S_1)$  is an $I$-complex.  Conversely, if $\Gamma$ is an $I$-complex over every field $k,$ then $\Gamma$ is such a complex.
\end{thm}

\begin{proof}
  Let $\Gamma$ be a bounded subcomplex of $L(\A,S_1)$ and let $\phi$ be a pointed lattice isomorphism from $(L(\A),S_1)$ to $(L(M^\star),e).$ For notational simplicity we assme that $E = \{1,\dots,n\} = [n]$ and $e=1.$ As usual, for each cell $c$ of $\Gamma$ let $\psi(c) = \max\{X \in \A: c \subseteq X\}.$ Since $\A$ is essential $\psi(c)$ is the meet of the coatoms of $L(\A)$ which correspond to the vertices of $c.$ If $v$ is a vertex in $\Gamma,$ then $\phi(v)$ is a coatom of $L(M^\star)$ which does not contain $e.$  Similarly, for each cell $c$ of $\Gamma,\ \phi(\psi(c))$ is a flat of $L(M^\star)$ which does not contain $e.$  Label each cell $c$ with $m_c,$ the square-free monomial whose support is $([n] - \phi(\psi(c ))) - \{1\}.$    Matroid duality implies that the support of $m_c$ is the union of the circuits of $M$ which are the complements of the coatoms corresponding to the vertices incident to $c.$ Thus each cell of $\Gamma$ is labeled with $\lcm\{m_i:v_i \le c\}.$  As the rank of any flat $X$ of $L(\A)$ is equal to $d-\dim(X), \Gamma$ satisfies $m_c \neq m_{c^\prime}$ whenever $c > c^\prime.$ Applying matroid duality again, we see that every $I$-essential subcomplex is of the form $\phi^{-1}(X) \cap \Gamma,$ where $X \in L(M^\star,e).$ The regularity of $\A$ with respect to $S_1$ guarantees that every $I$-essential subcomplex of $\Gamma$ is contractible, and hence acyclic.

For the converse, assume that $\Gamma$ is an $I$-complex over every field $k.$ Then each $I$-essential subcomplex of $\Gamma$ is acyclic over $\Z.$  Relabel each cell $c$  with the complement in $[n]-\{1\}$ of the support of $m_c$ i.e., the flat of $M^\star$ which does not contain $e$ and corresponds by matroid duality to $\lcm_{v \le c} m_v.$  The zero-skeleton of $\Gamma$ is the same as the zero-skeleton of any bounded subcomplex of $(\A,S_1)$ where $\A$ is an essential arrangement of homotopy spheres which is regular with respect to $S_1$ and $(L(\A),S_1) \cong (L(M^\star),e).$  As noted above, all $I$-essential subcomplexes of $\Gamma$ consist of cells whose labels contain a fixed $X \in L(M^\star,e).$ Proceeding inductively on the corank of all the flats in $L(M^\star,e)$ we see that $\Gamma$ must be equivalent to one constructed in exactly the same fashion as the procedure in the representation theorem for constructing $D_b$ for $L(M^\star,e).$  So, $D_b$ and the simultaneously constructed $(\A,S_1)$ are the required complexes.

\end{proof}

\begin{cor} \label{betas}
  Let $I = I_{\S(M,e)}$ be a matroid Steiner ideal.  Then
  $$\beta_i(I) = |\sum_{\stackrel{r(X) = n-r-i}{e \nleq Y}} \mu_{L(M^\star)}(X,Y)|$$
\end{cor}

\begin{proof} Theorem \ref{I-complexes} and  Corollary \ref{bounded f-vector} \end{proof}

When $\Delta=\Delta(M)$ is the independence complex of a matroid, Corollary \ref{betas} recovers Stanley's formula \cite{St1}
$$\beta_i(I_\Delta) = |\sum_{r(X) = n-r-i} \mu_{L(M^\star)}(X,\hat{1})|.$$

\noindent In this case $\Delta = \S(\tM,\tilde{e}),$ where $\tM$ is the 
free coextension of $M$ and $\tilde{e}$ is the extra point.  The poset of flats of $\tM^\star$ which do not lie above $\tilde{e}$ is $L(M^\star)$ with $\hat{1}$ removed.  So, $$ \mu_{L(M^\star)}(X,\hat{1}) = - \sum_{\stackrel{Y \in L(\tM^\star)}{\tilde{e} \nleq Y}} \mu_{L(\tM^\star)}(X,Y).$$
 
The method developed by Novik et al. to construct  minimal resolutions for face ideals of independence complexes of matroids  can be applied to matroid Steiner ideals.  Instead of applying \cite[Corollary 3.10]{NPS} to the order dual of $L(M^\star),$ we can use the order dual of $\hat{L}(M^\star,e) = L(M^\star,e) \cup \{\hat{1}\}.$  The only non-trivial point is that every upper interval of $\hat{L}(M^\star,e)$ is Cohen-Macaulay.   This follows from the fact that for any $X \in \hat{L}(M^\star,e)$ the upper interval $[X,\hat{1}]$ is isomorphic to $\hat{L}(M^\star/X,e)$ and every geometric semilattice is shellable and hence Cohen-Macaulay \cite{WW}.

Equivalence classes of CW-complexes which are $I$-complexes over every field can also be parameterized algebraically. 
Following the notation of \cite[pg.\ 299]{NPS}, view the complex $\mathcal{Z}(P)$ as a complex over $\Z.$  Use $P$ equal to the order dual of $\hat{L}(M^\star,e).$ By Theorem  \ref{I-complexes} equivalence classes of such $I$-complexes come from all the possible $D_b$ constructed in the representation theorem for $(L(M^\star),e).$  Every $D_b$ corresponds to choosing  bases for $H_\star(\V_X)$ for each $X \in L(M^\star,e).$  Working backwards from $Z_{-1}$ and using the fact that $\phi$ is an injection when restricted to each direct sum component, this is equivalent to choosing a $\Z$-basis for every direct sum component which occurs in $\mathcal{Z}(P).$  The only restriction to these bases is that the image under $\phi$ of any basis of a component in $\bigoplus_{rk(F) = 2} Z_0(\Delta(F))$ must be the difference of exactly two basis elements of $\bigoplus_{rk(F)=1}Z_{-1}(\Delta(F)).$  This is a reflection of the fact that the boundary of any one-cell is always the difference of two zero-cells. \\

\noindent {\it Acknowledgements:} Bernd Sturmfels suggested the problem of determining all $I$-complexes when $I$ is the face ideal of the independence complex of a matroid. Louis Billera and Tom Zaslavsky helped clear up several points in the exposition and the latter also pointed out the relevance of \cite{WW}.

\bibliography{matroids3,geometry} 
\bibliographystyle{plain}

 \end{document}